\documentclass[12pt]{article} 
\usepackage{amssymb,amsmath}
\usepackage{amsthm}

\newtheorem{theorem}{Theorem}[section]

\newtheorem{lemma}[theorem]{Lemma}

\theoremstyle{definition}

\theoremstyle{remark}

\begin{document}
\newcommand{\beq}{\begin{equation}} 
\newcommand{\eeq}{\end{equation}}
\newcommand{\zz}{\mathbb{Z}}
\newcommand{\pp}{\mathbb{P}} 
\newcommand{\nn}{\mathbb{N}}
\newcommand{\qq}{\mathbb{Q}}
\newcommand{\rr}{\mathbb{R}}
\newcommand{\parn}{\mathrm{Par}(n)}
\newcommand{\bm}[1]{{\mbox{\boldmath $#1$}}}
\newcommand{\con}{\mathrm{Comp}(n)}
\newcommand{\sn}{\mathfrak{S}_n} 
\newcommand{\fs}{\mathfrak{S}}
\newcommand{\st}{\,:\,} 
\newcommand{\as}{\mathrm{as}}
\newcommand{\is}{\mathrm{is}}
\newcommand{\lgn}{\mathrm{len}}
\newcommand{\dis}{\displaystyle}
\newcommand{\bea}{\begin{eqnarray}}
\newcommand{\eea}{\end{eqnarray}}
\newcommand{\be}{\begin{enumerate}}
\newcommand{\ee}{\end{enumerate}}
\newcommand{\beas}{\begin{eqnarray*}}
\newcommand{\eeas}{\end{eqnarray*}}
\newcommand{\ffac}[2]{\langle #1\rangle_{#2}}

\thispagestyle{empty}

\vskip 20pt
\begin{center}
{\large\bf Some Combinatorial Properties of Hook Lengths, Contents,
  and Parts of Partitions}$^1$
\vskip 15pt
{\bf Richard P. Stanley}$^2$\\[.1in]
{\it Department of Mathematics}\\ 
{\it Massachusetts Institute of Technology}\\
{\it Cambridge, MA 02139, USA}\\
{\texttt{rstan@math.mit.edu}}\\[.2in]
{\bf\small version of 25 March 2009}\\[.2in]
{\textsc{Dedicated to George Andrews for his 70th birthday}}
\end{center}

 \begin{abstract}
   The main result of this paper is a generalization of a conjecture
   of Guoniu Han, originally inspired by an identity of Nekrasov and
   Okounkov. Our result states that if $F$ is any symmetric function
   (say over $\qq$) and if
  $$ \Phi_n(F) = \frac{1}{n!}\sum_{\lambda\vdash n} f_\lambda^2\,
            F(h_u^2\st u\in\lambda), $$
where $h_u$ denotes the hook length of the square $u$ of the partition
$\lambda$ of $n$ and $f_\lambda$ is the number of standard Young
tableaux of shape $\lambda$, then $\Phi_n(F)$ is a polynomial function
of $n$. A similar result is obtained when $F(h_u^2\st u\in\lambda)$ is
replaced with a function that is symmetric separately in the contents
$c_u$ of $\lambda$ and the shifted parts $\lambda_i+n-i$ of $\lambda$.
 \end{abstract}

\section{Introduction.} \label{sec1} We assume basic knowledge of
symmetric functions such as given in \cite[Ch.~7]{ec2}. Let
$f_\lambda$ denote the number of standard Young tableaux (SYT) of
shape $\lambda\vdash n$. Recall the hook length formula of Frame,
Robinson, and Thrall \cite{f-r-t}\cite[Cor.~7.21.6]{ec2}:
  $$ f_\lambda = \frac{n!}{\prod_{u\in\lambda} h_u}, $$
where $u$ ranges over all squares in the (Young) diagram of $\lambda$,
and $h_u$ denotes the hook length at $u$. A basic property of the
numbers $f_\lambda$ is the formula
  $$ \sum_{\lambda\vdash n} f_\lambda^2 = n!, $$
which has an elegant bijective proof (the RSK algorithm). We will be
interested in generalizing this formula by weighting the sum on the
left by various functions of $\lambda$. Our primary interest is the sum
  $$ \Phi_n(F) = \frac{1}{n!}\sum_{\lambda\vdash n} f_\lambda^2\,
            F(h_u^2\st u\in\lambda), $$
where $F=F(x_1,x_2,\dots)$ is a symmetric function, say over $\qq$
(denoted $F\in \Lambda_\qq$).  The notation $F(h_u^2\st u\in\lambda)$
means that we are substituting for $n$ of the variables in $F$ the
quantities $h_u^2$ for $u\in\lambda$, and setting all other variables
equal to 0. For instance, if $F=p_k:= \sum x_i^k$, then
  $$ F(h_u^2\st u\in\lambda) = \sum_{u\in\lambda} h_u^{2k}. $$

This paper is motivated by the conjecture \cite[Conj.~3.1]{han} of
Guoniu Han that for all $k\in\pp=\{1,2,\dots\}$, we have that
$\Phi_n(p_k)\in\qq[n]$, i.e.,  
  $$ \frac{1}{n!}\sum_{\lambda\vdash n} f_\lambda^2\sum_{u\in
       \lambda} h_u^{2k} $$
is a polynomial function of $n$. This conjecture in turn was inspired
by the remarkable identity of Nekrasov and Okounkov \cite{n-o} (later
given a more elementary proof by Han \cite{han2})
  \beq \sum_{n\geq 0}\left(\sum_{\lambda\vdash n}f_\lambda^2 
   \prod_{u\in\lambda} (t+h_u^2)\right)\frac{x^n}{n!^2} =
   \prod_{i\geq 1}(1-x^i)^{-1-t}. \label{eq:n-o} \eeq
(We have stated this identity in a slightly different form than given
in \cite{han2}\cite{n-o}.) Our main result
(Theorem~\ref{thm:main}) states that $\Phi_n(F)\in\qq[n]$ for any
$F\in\Lambda_{\qq}$, i.e., for fixed $F$, $\Phi_n(F)$ is a polynomial
function of $n$. In the course of the proof we also show that 
  $$ \frac{1}{n!}\sum_{\lambda\vdash n}
    f_\lambda^2\,G(\{c_u\st u\in\lambda\}; \{ \lambda_i+n-i\st 1\leq
    i\leq n\})\in \qq[n]. $$
Here $G=G(x;y)$ is any formal power series of bounded degree over
$\qq$ that is symmetric in the $x$ and $y$ variables separately.
Moreover, $c_u$ denotes the content of $u\in\lambda$
\cite[p.~373]{ec2}; and we write
$\lambda=(\lambda_1,\dots,\lambda_n)$, adding 0's at the end so that
there are exactly $n$ parts.

\medskip 
\textsc{Acknowledgment.} I am grateful to Soichi Okada for calling my
attention to reference \cite{fkmo} and for providing
conjecture~\eqref{eq:okconj}. I also am grateful to an anonymous
referee for many helpful suggestions.

\section{Contents.} \label{sec2} 
In the next section we will obtain a stronger result than the main
result of this section (Theorem~\ref{thm:cu}). Since
Theorem~\ref{thm:cu} may be of independent interest and may be helpful
for understanding the next section, we treat it separately.

If $t\in\pp$ and $F$ is a symmetric function in the variables
$x_1,x_2,\dots$, then we write $F(1^t)$ for the result of setting
$x_1=x_2=\cdots=x_t=1$ and all other $x_j=0$ in $F$. For instance,
$p_\lambda(1^t)= t^{\ell(\lambda)}$, where $\ell(\lambda)$ is the
number of (positive) parts of $\lambda$. The \emph{hook-content}
formula for the case $q=1$ \cite[Cor.~7.21.4]{ec2} asserts that
  \beq s_\lambda(1^t) =\frac{\prod_{u\in\lambda}(t+c_u)}{H_\lambda}, 
    \label{eq:hoco} \eeq
where $s_\lambda$ is a Schur function and 
  $$ H_\lambda=\prod_{u\in\lambda} h_u, $$
the product of the hook lengths of $\lambda$ (so
$f_\lambda=n!/H_\lambda$). 

\begin{theorem} \label{thm:cu}
For any $F\in\Lambda_{\qq}$ we have 
  $$ \frac{1}{n!}\sum_{\lambda\vdash n} f_\lambda^2\,
            F(c_u\st u\in\lambda) \in \qq[n]. $$
\end{theorem}

\proof
By linearity it suffices to take $F=e_\mu$, the elementary symmetric
function indexed by $\mu$. Let $k\in\pp$, and for
$1\leq i\leq k$ let $x^{(i)}$ denote the set of 
variables $x_1^{(i)}, x_2^{(i)}, \dots$. Let $\sn$ denote the
symmetric group of all permutations of $\{1,\dots,n\}$. For $w\in\sn$
write $\rho(w)$ for the cycle type of $w$, i.e., $\rho(w)$ is the
partition of $n$ whose parts are the cycle lengths of $w$. We use the
identity \cite[Prop.~2.2]{h-s-s}\cite[Exer.~7.70]{ec2} 
  \beq \sum_{\lambda\vdash n} H_\lambda^{k-2} s_\lambda(x^{(1)})
     \cdots s_\lambda(x^{(k)}) = \frac{1}{n!}\sum_{\substack
      {w_1 w_2\cdots w_k=1\\ \mathrm{in}\ \sn}} p_{\rho(w_1)}
      (x^{(1)})\cdots  p_{\rho(w_k)}(x^{(k)}). \label{eq:hkm2} \eeq
Make the substitution $x^{(i)}=1^{t_i}$ as explained above. Letting
$c(w)$ denote the number of cycles of $w\in\sn$, we obtain
  \beq \sum_{\lambda\vdash n} H_\lambda^{-2}\prod_{u\in\lambda}
    \prod_{i=1}^k (t_i+c_u) = \frac{1}{n!} \sum_{\substack
      {w_1 w_2\cdots w_k=1\\ \mathrm{in}\ \sn}} t_1^{c(w_1)}
      \cdots t_k^{c(w_k)}. \label{eq:hcsub} \eeq
For any $n\geq \mu_1$ let $\mu=(\mu_1,\dots,\mu_k)$ be a partition
with $k$ parts, and take the coefficient of $t_1^{n-\mu_1} \cdots
t_k^{n-\mu_k}$ on both sides of equation~\eqref{eq:hcsub}. Using
$f_\lambda = n!/H_\lambda$, we obtain
  $$ \hspace{-2in}\frac{1}{n!}\sum_{\lambda\vdash n} f_\lambda^2\,
  e_\mu(c_u\st u\in\lambda) $$
 \beq \qquad = \#\{(w_1,\cdots,w_k)\in\sn^k\st w_1\cdots w_k=1,\ 
    c(w_i)=n-\mu_i\}. \label{eq:conid} \eeq
We therefore need to show that the right-hand side of
equation~\eqref{eq:conid} is a polynomial function of $n$.

Suppose that $c(w_i)=n-\mu_i$ and that the union $F$ of the non-fixed
points of all the $w_i$'s has $r$ elements. Then
   \beq 1+\mu_1 \leq r\leq 2\sum\mu_i. \label{eq:mui} \eeq 
We can choose the set $F$ in $\binom{n}{r}$ ways.  Once we make this
choice there is a certain number of ways (depending on $r$ but
independent of $n$) that we can have $w_1\cdots w_k=1$. (In more
algebraic terms, $\sn$ acts on $S_\mu$ by conjugation, where $S_\mu$
is the set on the right-hand side of \eqref{eq:conid}, and the number
of orbits of this action is independent of $n$.) Hence for $n\geq
1+ \mu_1$, $\#S_\mu$ is a finite linear combination (over
$\nn=\{0,1,2,\dots\}$) of polynomials $\binom nr$, and is thus a
polynomial $N_\mu(n)$ as desired.

If $n< 1+\mu_1$, then it is clear from the previous paragraph
that the polynomial $N_\mu$ satisfies $N_\mu(n)=0$. On the other hand,
if $\lambda\vdash n$ then we also have $e_\mu(c_u\st u\in\lambda)=0$.
Hence the two sides of equation~\eqref{eq:conid} agree for $0\leq n<
1+\max \mu_i$, and the proof is complete.  \qed

\medskip
Note that the proof of Theorem~\ref{thm:cu} shows that $N_\mu(n)$
is a \emph{nonnegative} integer linear combination of the polynomials
$\binom nr$. It can be shown that either $N_\mu=0$ or $\deg N_\mu=\sum
\mu_i$. Moreover $N_\mu\neq 0$ if and only $\sum \mu_i$ is even, say
$2r$, and $\mu_1\leq r$. The nonzero polynomials $N_\mu(n)$ for
$|\mu|\leq 6$ are given by  
  \beas N_{1,1}(n) & = & \frac{n(n-1)}{2}\\[1em]
     N_{2,2}(n) & = & \frac{n(n-1)(n-2)(3n-1)}{24}\\[1em]
     N_{2,1,1}(n) & = & \frac{n(n-1)(n-2)(n+1)}{4}\\[1em]
     N_{1,1,1,1}(n) & = & \frac{n(n-1)(3n^2+n-12)}{4}\\[1em]
     N_{3,3}(n) & = & \frac{n^2(n-1)^2(n-2)(n-3)}{48}\\[1em]
     N_{3,2,1}(n) & = &
        \frac{n(n-1)(n-2)(n-3)(3n^2+5n+4)}{48}\\[1em] 
     N_{3,1,1,1}(n) & = & \frac{n(n-1)(n-2)(n-3)(n^2+3n+4)}{8}\\[1em]
     N_{2,2,2}(n) & = & \frac{n(n-1)(n-2)(3n^3-9n-46)}{24}\\[1em]
     N_{2,2,1,1}(n) & = & \frac{n(n-1)(n-2)(15n^3+20n^2-59n-312)}
        {48} \eeas
   \beas
     N_{2,1,1,1,1}(n) & = & \frac{n(n-1)(n-2)(3n^3+8n^2-7n-96)}{4}\\[1em]
     N_{1,1,1,1,1,1}(n) & = & \frac{n(n-1)(15n^4+30n^3-105n^2-700n+1344)}
       {8}. \eeas

A slight modification of the proof of a special case of
Theorem~\ref{thm:cu} leads to a ``content Nekrasov-Okounkov formula.''

\begin{theorem} \label{thm:cno}
We have
  $$ \sum_{n\geq 0}\left( \sum_{\lambda\vdash n} f_\lambda^2
  \prod_{u\in\lambda} (t+c_u^2)\right)\frac{x^n}{n!^2} =
   (1-x)^{-t}. $$ 
\end{theorem} 

\proof
By the ``dual Cauchy identity'' \cite[Thm.~7.14.3]{ec2} we have
  $$ \sum_{\lambda\vdash n} s_\lambda(x)s_{\lambda'}(y) = \frac{1}{n!}
    \sum_{w\in\sn} \varepsilon_w p_{\rho(w)}(x)p_{\rho(w)}(y), $$
where $\varepsilon(w)$ is given by
equation~\eqref{eq:varep}, and where $\lambda'$ denotes the conjugate
partition to $\lambda$. Substitute $x=1^t$ and $y=1^t$. Since the
contents of $\lambda'$ are the negative of those of $\lambda$, we
obtain
 $$ \sum_{\lambda\vdash n} H_\lambda^{-2}\prod_{u\in\lambda}
   (t^2-c_u^2) = \frac{1}{n!}\sum_{w\in\sn}\varepsilon_w t^{2c(w)}. $$
It is a well-known and basic fact that the sum on the right is
$\binom{t^2}{n}$. Put $-t$ for $t^2$, multiply by $(-x)^n$ and
sum on $n\geq 0$ to get the stated formula.
\qed

\medskip
A simple variant of Theorem~\ref{thm:cno} follows from considering the
usual Cauchy identity (the case $k=2$ of equation~\eqref{eq:hkm2})
instead of the dual one:
  $$  \sum_{n\geq 0}\left( \sum_{\lambda\vdash n} f_\lambda^2
  \prod_{u\in\lambda} (t+c_u)(v+c_u)\right)\frac{x^n}{n!^2} =
   (1-x)^{-tv}. $$ 

A related identity is due to Fujii \emph{et al.}\
\cite[Appendix]{fkmo}, namely, for any $r\geq 0$ we have
 \beq \frac{1}{n!} \sum_{\lambda \vdash n}
  ( f^\lambda )^2 \sum_{u \in \lambda}
    \prod_{i=0}^{r-1} ( c_u^2 - i^2 ) =
    \frac{(2r)!}{(r+1)!^2} \ffac{n}{r+1}, \label{eq:okada} \eeq
where $\ffac{n}{r+1}=n(n-1)\cdots (n-r)$. It follows from this formula 
that
  \beq \frac{1}{n!} \sum_{\lambda \vdash n}
  ( f^\lambda )^2 \sum_{u \in \lambda}c_u^{2k} =
   \sum_{j=1}^k T(k,j)\frac{(2j)!}{(j+1)!^2}\ffac{n}{j+1}, 
    \label{eq:okada2} \eeq
where $T(k,j)$ is a \emph{central factorial number}
\cite[Exer.~5.8]{ec2}. One of several equivalent definitions of
$T(k,j)$ is the explicit formula
  $$ T(k,j) =2\sum_{i=1}^j \frac{(-1)^{j-i}i^{2k}}{(j-i)!(j+i)!}. $$ 
Another definition is the generating function 
  \beq \sum_{k\geq 0} T(k,j)x^k = \frac{x^j}{(1-1^2x)(1-2^2x)
     \cdots(1-j^2x)}. \label{eq:tgf} \eeq
The equivalence of equations~\eqref{eq:okada} and \eqref{eq:okada2} is
a simple consequence of \eqref{eq:tgf}.
For ``hook length analogues'' of
equations~\eqref{eq:okada} and \eqref{eq:okada2}, see the Note at the
end of Section~\ref{sec4}. 

\section{Shifted parts.} \label{sec3}
In this section we write partitions $\lambda$ of $n$ as
$(\lambda_1,\dots,\lambda_n)$, placing as many 0's at the end as
necessary. Thus for instance the three partitions of 3 are $(3,0,0)$,
$(2,1,0)$, and $(1,1,1)$. Let $G(x;y)$ be a formal power
series over $\qq$ of bounded degree that is
symmetric in the variables $x=(x_1,x_2,\dots)$ and $y=(y_1,y_2,\dots)$
separately; in symbols, $G\in \Lambda_\qq[x]\otimes
\Lambda_\qq[y]$. We are interested in the quantity 
  \beq \Psi_n(G) = \frac{1}{n!}\sum_{\lambda\vdash n}
    f_\lambda^2\,G(\{c_u\st u\in\lambda\}; \{ \lambda_i+n-i\st 1\leq
    i\leq n\}). \label{eq:psing} \eeq
       The case $y_i=0$ for all $i$ reduces to what was considered in the
previous section.  We will show that $\Psi_n(G)$ is a polynomial
in $n$ by an argument similar to the proof of Theorem~\ref{thm:cu}. In
addition to the substitution $x^{(i)}=1^{t_i}$ we use a certain linear
transformation $\varphi$ which we now define.

Let $x^{(1)}, \dots,x^{(j)}$ and $y^{(1)}, \dots, y^{(k)}$ be disjoint
sets of variables. We will work in the ring $R$ of all bounded formal
power series over $\qq$ that are symmetric in each set of variables
separately. Define a map $\varphi\colon R\rightarrow
\qq[v_1,\dots,v_k]$ by the conditions:
 \begin{itemize}
  \item The map $\varphi$ is linear over $\Lambda_\qq[x^{(1)}] \otimes
 \cdots\otimes
 \Lambda_\qq[x^{(j)}]$, i.e, the $x^{(i)}$-variables are treated as
  scalars.
  \item We have 
   $$ \varphi\left(s_\lambda(y^{(h)})\right) =
    \frac{\prod_{i=1}^n (v_h+\lambda_i+n-i)}{H_\lambda}, $$
 where $\lambda\vdash n$.
  \item We have
   $$   \varphi\left(G_1(y^{(1)})\cdots G_k(y^{(k)})\right) =
    \varphi\left(G_1(y^{(1)})\right)\cdots 
      \varphi\left(G_k(y^{(k)})\right), $$
where $G_h\in\Lambda_{\qq}[x^{(1)},\dots,x^{(j)}, y^{(h)}]$.
  \end{itemize}    
More algebraically, let $\Psi= \Lambda_\qq[x^{(1)}] \otimes
 \cdots\otimes \Lambda_\qq[x^{(j)}]$, and let $\varphi_h\colon
 \Psi[y^{(h)}]\rightarrow \qq[v_h]$ be the $\Psi$-linear
 transformation defined by 
  $$ \varphi_h(s_\lambda(y^{(h)})) = H_\lambda^{-1}\prod_{i=1}^n
   (v_h+\lambda_i+n-i). $$ 
Then $\varphi = \varphi_1\otimes\cdots
\otimes \varphi_k$ (tensor product over $\Psi$).

Write for simplicity $f$ for $f(y^{(1)})$ and $v$ for $v_1$.
We would like to evaluate $\varphi(p_\mu)$, where $p_\mu$ is a
power-sum symmetric function. We first need the following
lemma. Define
  $$ A_\lambda(v)=H_\lambda^{-1} (v+\lambda_1+n-1)
       (v+\lambda_2+n-2)\cdots(v+\lambda_n). $$

\begin{lemma} \label{lemma:spid}
For all $n\geq 0$ we have
 \beq   \sum_{i=0}^n \binom{v+i-1}{i} p_1^i e_{n-i} =
  \sum_{\lambda\vdash n} A_\lambda(v) s_\lambda. 
    \label{eq:pelem} \eeq
Equivalently, we have
  $$ (1-p_1)^{-v}\sum_{n\geq 0}e_n = \sum_{n\geq 0}
    \sum_{\lambda\vdash n} A_\lambda(v) s_\lambda. $$
\end{lemma}
 
\emph{First proof} (sketch). I am grateful to Guoniu Han for providing
the following proof. Complete details may be found in his paper
\cite{han3}.  Denote the left-hand side of equation~\eqref{eq:pelem}
by $L_n(v)$ and the right-hand side by $R_n(v)$. It is easy to see that 
$L_n(v)=L_n(v-1)+p_1L_{n-1}(v)$, $L_n(0)=R_n(0)$, and
$L_0(v)=R_0(v)$. Hence we need to show that 
  \beq R_n(v)=R_n(v-1)+p_1R_{n-1}(v). \label{eq:rnrec} \eeq
Now for $\lambda\vdash n$ let
  $$ E_\lambda(v) =  A_\lambda(v+n+1)-A_\lambda(v+n)
  -\sum_{\mu\in\lambda\backslash 1}A_\mu(v+n+1), $$
where $\lambda\backslash 1$ denotes the set of all partitions $\mu$
obtained from $\lambda$ by removing one corner. Clearly $E_\lambda(v)$
is a polynomial in $v$ of degree at most $n$, and it is not difficult
to check that the degree in fact is at most $n-2$. The core of the
proof (which we omit) is to show that $E_\lambda(i-\lambda_i)=0$ for
$i=1,2,\dots,n-1$. Since $E_\lambda(v)$ has degree at most $n-2$ and
vanishes at $n-1$ distinct integers, we conclude that
$E_\lambda(v)=0$. It is now straightforward to verify that
equation~\eqref{eq:rnrec} holds.
\qed

\medskip
\emph{Second proof.} I am grateful to Tewodros Amdeberhan for helpful
discussions. A formula of Andrews, Goulden, and Jackson
\cite{a-g-j} asserts that
  $$ \hspace{-.5in}\sum_\lambda s_\lambda(y_1,\dots,y_n)
     s_\lambda(z_1,\dots,z_m) \prod_{i=1}^n(v-\lambda_i-n+i) $$
  $$ = \prod_{j=1}^n\prod_{k=1}^m \frac{1}{1-y_jz_k}\cdot[t_1\cdots
  t_n] (1+t_1+\cdots+t_n)^v \prod_{k=1}^m\left( 1-\sum_{j=1}^n
       \frac{t_j y_jz_k}{1-y_jz_k}\right), $$
where the sum is over all partitions $\lambda$ satisfying
$\ell(\lambda)\leq n$, and where $[t_1\cdots t_n]X$ denotes the
coefficient of $t_1\cdots t_n$ in $X$. Change $v$ to $-v$ and multiply
by $(-1)^n$ to get 
  $$ \hspace{-.5in}\sum_\lambda s_\lambda(y_1,\dots,y_n)
     s_\lambda(z_1,\dots,z_m) \prod_{i=1}^n(v+\lambda_i+n-i) $$
  $$ \hspace{-10em}= (-1)^n\prod_{j=1}^n\prod_{k=1}^m
        \frac{1}{1-y_jz_k}\cdot $$ 
  $$ \qquad\quad  [t_1\cdots
     t_n] (1+t_1+\cdots+t_n)^{-v} \prod_{k=1}^m\left( 1-\sum_{j=1}^n
       \frac{t_j y_jz_k}{1-y_jz_k}\right). $$
Let $m=n$, and take the
coefficient of $z_1\cdots z_n$ on both sides. The
left-hand side becomes
  $$ \sum_{\lambda\vdash n} f_\lambda s_\lambda(y)\prod_{i=1}^n
    (v+\lambda_i+n-i). $$
 Consider the coefficient of $z_1\cdots z_n$ on the
 right-hand side. A term from this coefficient is obtained as
 follows. Pick a subset $S$ of $[n]=\{1,2,\dots,n\}$, say
 $\#S=r$. Choose the 
 coefficient of $\prod_{i\in S}z_i$ from $\prod_{j=1}^n\prod_{k=1}^n
 (1-y_jz_k)^{-1}$. This coefficient is $p_1(y)^r$, and there are
 $\binom nr$ choices for $S$. We now must choose the coefficient 
 $\prod_{i\in [n]-S}z_i$ from $\prod_{k=1}^n\left( 1-\sum_{j=1}^n
 \frac{t_j y_jz_k}{1-y_jz_k}\right)$. This coefficient is $(-1)^{n-r}
 (t_1y_1+\cdots+t_ny_n)^{n-r}$. Hence 
  $$ \sum_{\lambda\vdash n} f_\lambda s_\lambda(y)\prod_{i=1}^n
    (t+\lambda_i+n-i) $$
  $$ = (-1)^n \sum_{r=0}^n \binom nr p_1(y)^r [t_1\cdots t_n]
     \frac{(-1)^{n-r}
       (t_1y_1+\cdots+t_ny_n)^{n-r}}{(1+t_1+\cdots+t_n)^{-v}}. $$  
Let $\{i_1,\dots,i_{n-r}\}$ be an $(n-r)$-element subset of $[n]$, and
let $\{j_1,\dots,j_r\}$ be its complement. Then
  \beas [t_{i_1} \cdots t_{i_{n-r}}]
  (t_1y_1+\cdots+t_ny_n)^{n-r} & = & (n-r)!y_{i_1}\cdots
  y_{i_{n-r}}\\[1em] 
    [t_{j_1}\cdots t_{j_r}] (1+t_1+\cdots+t_n)^{-v} & = &
    \binom{-v}{r}r!. \eeas
Hence
  $$ \hspace{-2in} \sum_{\lambda\vdash n} f_\lambda
  s_\lambda(y)\prod_{i=1}^n (t+\lambda_i+n-i) $$
  \beq \qquad = \sum_{r=0}^n r!(n-r)!\binom nr p_1(y)^r
      (-1)^r\binom{-v}{r} e_{n-r}(y). \label{eq:spb} \eeq
Write $(-1)^r\binom{-v}{r} = \binom{v+r-1}{r}$ and divide both sides
of equation~\eqref{eq:spb} by $n!$ to complete the proof.
\qed

\medskip
\textsc{Note.} (a) Amdeberhan \cite{tew} has simplified the second
proof of Lemma~\ref{lemma:spid}; in particular, he avoids the use of
the Andrews-Goulden-Jackson formula.

(b) Since the left-hand side of equation~\eqref{eq:pelem}
is an \emph{integral} linear combination of Schur functions when
$v\in\zz$ (e.g., by Pieri's rule), it follows that for every 
$v\in\zz$ we have $A_\lambda(v)\in\zz$.
By expanding the left-hand side of \eqref{eq:pelem} in terms of Schur
functions, we in fact obtain the following combinatorial expression
for $A_\lambda(v)$:
  $$ A_\lambda(v) = \sum_{i=0}^n\binom{v+i-1}{i}f_{\lambda/1^{n-i}},
  $$
where $f_{\lambda/1^{n-i}}$ denotes the number of SYT of the skew
shape $\lambda/1^{n-i}$.

\medskip
We now turn to the evaluation of $\varphi(p_\mu)$. 

\begin{lemma} \label{lemma:vphi}
For any partition $\mu\vdash n$ with $\ell=\ell(\mu)$ nonzero parts,
we have  
  $$ \varphi(p_\mu) = (-1)^{n-\ell} \sum_{i=0}^m \binom mi
    (v)_i, $$
where $m=m_1(\mu)$, the number of parts of $\mu$ equal to 1, and
$(v)_i=v(v+1)\cdots(v+i-1)$.
\end{lemma}
     
\proof
We will work with two sets of variables $x=(x_1,x_2,\dots)$ and
$y=(y_1, y_2,\dots)$. Recall that $\varphi$ acts on symmetric
functions in $y$ only, regarding symmetric function in $x$ as
scalars. Thus using Lemma~\ref{lemma:spid} we have
  \bea \varphi \sum_{\lambda\vdash n}s_\lambda(x)s_\lambda(y) & = &
     \sum_{\lambda\vdash n} A_\lambda(v) s_\lambda(x). \nonumber\\ & = &
     \sum_{i=0}^n \binom{v+i-1}{i}p_1^ie_{n-i}. \label{eq:varphi} \eea
A standard symmetric function identity \cite[(7.23)]{ec2} states that
  $$ e_{n-i} =\sum_{\rho\vdash n-i}\varepsilon_\rho z_\rho^{-1}p_\rho,
  $$
where 
  \beq \varepsilon_\rho=(-1)^{|\rho|-\ell(\rho)}, \label{eq:varep}
  \eeq 
and if $\rho$ has $m_i$ parts equal to $i$ then
$z_\rho=1^{m_1}m_1!2^{m_2}m_2!\cdots$.  Let $\nu$ be the partition
obtained from $\mu$ by removing all parts equal to 1. Write
$(\nu,1^j)$ for the partition obtained from $\nu$ by adjoining $j$
1's, so $\mu=(\nu,1^m)$. Note that
 $$ \varepsilon_{(\nu,1^{m-i})}  = 
 (-1)^{|\nu|+m-i-\ell(\nu)-(m-i)} = (-1)^{|\nu|-\ell(\nu)}
    = (-1)^{n-\ell(\mu)}. $$
Note also that 
  $$ z_{(\nu,1^{m-i})}=\frac{(m-i)!}{m!}z_\mu. $$
Hence if we expand the right-hand side of equation~\eqref{eq:varphi}
in terms of power sum symmetric functions, then the coefficient of
$p_\mu$ is
  $$ \sum_{i=0}^m \binom{v+i-1}{i}\varepsilon_{(\nu,1^{m-i})} 
         z_{(\nu,1^{m-i})}^{-1} $$
   \beq \qquad  = (-1)^{n-\ell}
    \sum_{i=0}^m \binom mi (v)_i z_\mu^{-1}. 
     \label{eq:pcoef} \eeq

It follows from the Cauchy identity \cite[Thm.~7.12.1]{ec2} (and is
also the special case $k=2$ of equation~\eqref{eq:hkm2}) that
  \beq \sum_{\lambda\vdash n} s_\lambda(x)s_\lambda(y) =
     \sum_{\mu\vdash n} z_\mu^{-1} p_\mu(x)p_\mu(y). 
    \label{eq:cauchy} \eeq
Thus when we apply $\varphi$ (acting on the $y$ variables) to
equation~\eqref{eq:cauchy} and use \eqref{eq:pcoef}, then we obtain
  $$ \hspace{-1in} \sum_{\mu\vdash n} \varphi(p_\mu(y))p_\mu(x) $$ 
   $$ \qquad =
    \sum_{\mu\vdash n} \left((-1)^{n-\ell(\mu)}
    \sum_{i=0}^m \binom mi (v)_i\right) z_\mu^{-1}
     p_\mu(x). $$
Since the $p_\mu$'s are linearly independent, the proof follows.
\qed

\begin{theorem} \label{thm:sp}
For any $G\in\Lambda_\qq[x]\otimes \Lambda_\qq[y]$ we have 
  $$ \Psi_n(G) \in \qq[n], $$
where $\Psi_n(G)$ is given by equation~\eqref{eq:psing}.
\end{theorem}

\proof
By linearity it suffices to take $G=e_\mu(x)e_\nu(y)$. Apply $\varphi$
to the identity \eqref{eq:hkm2} in the variables
$x^{(1)},\dots,x^{(j)}$, $y^{(1)},\dots,y^{(k)}$. Then make the
substitution $x^{(h)}=1^{t_h}$ and multiply by
$n!$. By equation~\eqref{eq:hoco} and Lemma~\ref{lemma:vphi} we obtain
  $$ \hspace{-.5in} \frac{1}{n!} \sum_{\lambda\vdash n} f_\lambda^2
      \prod_{h=1}^j \prod_{u\in\lambda} (t_h+c_u)\cdot
   \prod_{h=1}^k \prod_{i=1}^n (v_h+\lambda_i+n-i) $$
   $$ \qquad = 
              \sum_{\substack{w_1\cdots w_j w'_1\cdots w'_k=1\\
           \mathrm{in}\ \sn}} \prod_{h=1}^j t_h^{c(w_h)} $$
   \beq \cdot \prod_{h=1}^k \left( (-1)^{n-\ell(\rho(w'_h))}
        \sum_{i=0}^{m_1(\rho(w'_h))} \binom{m_1(\rho(w'_h))}{i}
        (v_h)_i\right). \label{eq:big} \eeq
The remainder of the proof is a straightforward generalization of that
of Theorem~\ref{thm:cu}. Take the coefficient of $t_1^{n-\mu_1}\cdots
t_j^{n-\mu_j}v_1^{n-\nu_1}\cdots v_k^{n-\nu_k}$. The left-hand side
becomes $\Psi_n(e_\mu(x)e_\nu(y))$,
so we need to show that the coefficient of $t_1^{n-\mu_1}\cdots
t_j^{n-\mu_j}v_1^{n-\nu_1}\cdots v_k^{n-\nu_k}$ on the right-hand side
of equation~\eqref{eq:big} is a polynomial in $n$.  Suppose that
$n\geq \mu_1$ and $n\geq\nu_1$. The coefficient of $v_h^{n-\nu_h}$ in
$v_h(v_h+1)\cdots (v_h+n-i-1)$ is the signless Stirling number
$c(n-i,n-\nu_h)$. The coefficient of $v_h^{n-\nu_h}$ in \eqref{eq:big}
is 0 unless $n-m_1(\rho(w'_h))\leq 
i\leq \nu_h$. For each choice of $0\leq i_h\leq i$ ($1\leq h\leq k$),
there are only finitely many orbits of the action of $\sn$ by
(coordinatewise) conjugation on the set of
$(w_1,\dots,w_j,w'_1,\cdots,w'_k)\in\sn^{j+k}$ for which $w_1\cdots
w_j w'_1\cdots w'_k=1$, $w_h$ has $n-\mu_h$ cycles, and $w'_h$ has
$n-i_h$ fixed points. The size of each of these orbits is a
polynomial in $n$, as in the proof of Theorem~\ref{thm:cu}. Moreover,
the Stirling number $c(n-i,n-\nu_h)$ is a polynomial in $n$ for fixed
$i$ and $\nu_h$, and similarly for the binomial coefficient
$\binom{n-i_h}{n-i}$, so $\Psi_n(e_\mu(x)e_\nu(y))$ is a
polynomial $N_{\mu,\nu}(n)$ for $n\geq \max\{\mu_1,\nu_1\}$.  If
$0\leq n<\max\{\mu_1,\nu_1\}$, then both $N_{\mu,\nu}(n)$ and
$\Psi_n(e_\mu(x)e_\nu(y))$ are equal to 0 (as in the proof of
Theorem~\ref{thm:cu}), so the proof is complete.  \qed

\medskip
\textsc{Note.} Since $n$ is a polynomial in $n$, it is easy to see
that Theorem~\ref{thm:sp} still holds if we replace $\Psi_n(G)$ with
  $$ \frac{1}{n!}\sum_{\lambda\vdash n}
    f_\lambda^2\,G(\{c_u\st u\in\lambda\}; \{ \lambda_i-i\st 1\leq
    i\leq n\}). $$
On the other hand, Theorem~\ref{thm:sp} becomes \emph{false} if we
replace $\Psi_n(G)$ with
    $$ \frac{1}{n!}\sum_{\lambda\vdash n} f_\lambda^2\, G(\{c_u\st
      u\in\lambda\}; \{\lambda_i\st 1\leq i\leq n\}). $$
For instance,
    $$ \frac{1}{n!}\sum_{\lambda\vdash n}
  f_\lambda^2(\lambda_1^2+\lambda_2^2+\cdots+\lambda_n^2) $$
is not a polynomial function of $n$, nor is it integer valued.

\section{Hook lengths squared.} \label{sec4}
The connection between contents, hook lengths, and the shifted parts
$\lambda_i+n-i$ is given by the following result, an
immediate consequence \cite[Lemma~7.21.1]{ec2}.

\begin{lemma} \label{lemma:mset}
Let $\lambda=(\lambda_1,\dots,\lambda_n)\vdash n$. Then we have the
multiset equality
  $$ \hspace{-1in} \{ h_u\st u\in \lambda\} \cup
  \{\lambda_i-\lambda_j-i+j\st 1\leq 
       i<j\leq n\} $$
  $$ \qquad = \{n+c_u\st u\in\lambda\}\cup
         \{1^{n-1},2^{n-2},\dots,n-1\}. $$
\end{lemma}

For example, when $\lambda=(3,1)$ Lemma~\ref{lemma:mset} asserts that
  $$ \{4,2,1,1\}\cup\{3,5,6,2,3,1\}=\{3,4,5,6\}\cup\{1,1,1,2,2,3\} $$
as multisets.

\begin{lemma} \label{lemma:mupo}
For any $F\in\Lambda_\qq$, we have
  $$ F(1^{n-1},2^{n-2},\dots,n-1)\in \qq[n], $$
where the exponents denote multiplicity.
\end{lemma}

\proof It suffices to take $F=p_j$ since the polynomials in $n$ form a
ring. Thus we want to show that
  $$ \sum_{i=1}^{n-1} (n-i)i^j \in\qq[n], $$
which is routine.
\qed

\medskip
We come to the main result of this paper. Recall the definition
  $$ \Phi_n(F) = \frac{1}{n!}\sum_{\lambda\vdash n} f_\lambda^2\,
            F(h_u^2\st u\in\lambda). $$

\begin{theorem} \label{thm:main}
For any symmetric function $F\in\Lambda_\qq$ we have $\Phi_n(F)
\in\qq[n]$. 
\end{theorem}

\proof
As usual it suffices to take $F=e_\mu$, where $\mu=(\mu_1,\dots,
\mu_k)$. Define the multisets (or \emph{alphabets})  
  \beas A_\lambda & = & \{h_u^2\st u\in\lambda\}\\
   B_\lambda & = & \{(\lambda_i-\lambda_j-i+j)^2\st 1\leq i<j\leq n\}\\
   C_\lambda & = & \{(n+c_u)^2\st u\in\lambda\}\\
   D_n & = & \{b_1^{n-1},b_2^{n-2},\dots, b_{n-1}\}, \eeas
where $b_i=i^2\in\zz$ (so for instance $D_4=\{1,1,1,4,4,9\}$). 
Write $\Omega(a,b,c)=(-1)^c e_a(C_\lambda)e_b(D_n)h_c(B_\lambda)$.
Using standard $\lambda$-ring notation and manipulations (see e.g.\
Lascoux \cite[Ch.~2]{lascoux}), we have from Lemma~\ref{lemma:mset}
that 
  \beas \Phi_n(e_\mu) & = & \frac{1}{n!}\sum_{\lambda\vdash n}
     f_\lambda^2\, e_\mu(A_\lambda)\\ & = &
     \frac{1}{n!}\sum_{\lambda\vdash n} 
     f_\lambda^2\, e_\mu(C_\lambda+D_n-B_\lambda)\\ & = & 
    \frac{1}{n!}\sum_{\lambda\vdash n} f_\lambda^2\,
    \prod_{i=1}^k \left( \sum_{\substack{a,b,c\geq 0\\ a+b+c=\mu_i}}
      \Omega(a,b,c)\right)\\ & = &
    \sum_{\substack{a_1,b_1,c_1\geq 0\\ a_1+b_1+c_1=\mu_1}}\cdots
     \sum_{\substack{a_k,b_k,c_k\geq 0\\ a_k+b_k+c_k=\mu_k}}
     \frac{1}{n!}\sum_{\lambda\vdash n} f_\lambda^2
     \prod_{r=1}^k \Omega(a,b,c). \eeas
Consider the inner sum over $\lambda$, together with the factor
$1/n!$.  By Lemma~\ref{lemma:mupo} each $e_{b_r}(D_n)$ is a polynomial
in $n$ which we can factor out of the sum. Note that
$h_{c_r}(B_\lambda)$ is a symmetric function of the numbers
$\rho_i=\lambda_i+n-i$ since $(\rho_i-\rho_j)^2$ is symmetric in $i$
and $j$. (This is the one point in the proof that requires the use of
the alphabet $\{h_u^2\st u\in\lambda\}$ rather than the more general
$\{h_u\st u\in\lambda$\}.)  What remains after factoring out each
$e_{b_r}(D_n)$ is therefore a polynomial in $n$ by
Theorem~\ref{thm:sp}, and the proof follows.  
\qed

\medskip
\textsc{Note.} (a) The $\lambda$-ring computations in the proof of
Theorem~\ref{thm:main} can easily be replaced with more ``naive''
techniques such as generating functions. The $\lambda$-ring approach,
however, makes the computation more routine.

(b) An interesting feature of the proofs of Theorems~\ref{thm:cu},
\ref{thm:sp}, and \ref{thm:main} is that they don't involve just
``formal'' properties of symmetric functions; use of representation
theory is required. This is because the only known proof of the
crucial equation~\eqref{eq:hkm2} involves representation theory, viz.,
the determination of the primitive orthogonal idempotents in the
center of the group algebra of $\sn$. Is there a proof of
\eqref{eq:hkm2} or of Theorems~\ref{thm:cu}, \ref{thm:sp}, and
\ref{thm:main} that doesn't involve representation theory?

\medskip
Here is a small table of the polynomials $\Phi_n(e_\mu)$:
  \beas \Phi_n(e_1) & = & \frac 12 n(3n-1)\\
   \Phi_n(e_2) & = & \frac{1}{24}n(n-1)(27n^2-67n+74)\\
   \Phi_n(e_1^2) &= & \frac{1}{12}n(27n^3-14n^2-9n+8)\\
   \Phi_n(e_3) & = & \frac{1}{48}n(n-1)(n-2)(27n^3-174n^2+511n-552)\\
   \Phi_n(e_2e_1) & = &
        \frac{1}{48}n(n-1)(81n^4-204n^3+137n^2+390n-512)\\
   \Phi_n(e_1^3) & = & \frac{1}{24}n(81n^5-45n^4-69n^3-31n^2+216n-128).
   \eeas
\indent\textsc{Note.} Soichi Okada has conjectured \cite{okada} the
    following ``hook analogue'' of equation~\eqref{eq:okada}:
 \beq
 \frac{1}{n!}
 \sum_{\lambda \vdash n}
  f_\lambda^2
  \sum_{u \in \lambda}
    \prod_{i=1}^r ( h_u^2 - i^2 ) =
 \frac{1}{2 (r+1)^2 } \binom{2r}{r} \binom{2r+2}{r+1} \ffac{n}{r+1}.
 \label{eq:okconj} \eeq
This conjecture has been proved by Greta Panova \cite{panova} using
Theorem~\ref{thm:main}. From this result we get the following
analogue of equation~\eqref{eq:okada2}:
  $$ \frac{1}{n!}\sum_{\lambda\vdash n} f_\lambda^2\sum_{u\in
       \lambda} h_u^{2k} = \sum_{j=1}^{k+1} T(k+1,j)\frac{1}{2j^2}
      \binom{2(j-1)}{j-1}\binom{2j}{j}\ffac{n}{j}. $$

\textsc{Note.} Using Theorem~\ref{thm:sp} and the method of the proof
of Theorem~\ref{thm:main} to reduce hook lengths squared to contents
  and shifted parts, it is clear that we have the following ``master
  theorem'' subsuming both Theorems~\ref{thm:sp} and \ref{thm:main}.

\begin{theorem} \label{thm:master}
For any $K\in\Lambda_\qq[x]\otimes\Lambda_\qq[y]\otimes\Lambda_q[z]$,
we have
  $$ \frac{1}{n!}\sum_{\lambda\vdash n}f_\lambda^2 K_\lambda \in
  \qq[n], $$
where 
 $$ K_\lambda=K(\{c_u\st u\in\lambda\}; \{ \lambda_i+n-i\st 1\leq
    i\leq n\}; \{ h_u^2\st u\in\lambda\}). $$ 
\end{theorem}

\section{Some questions.} \label{sec5}

\be
 \item  Can the Nekrasov-Okounkov formula \eqref{eq:n-o} be
proved using the techniques we have used to prove
Theorem~\ref{thm:main}? 
 \item Can the Nekrasov-Okounkov formula \eqref{eq:n-o} be generalized
with the left-hand side replaced with the following expression (or
some simple modification thereof)? 
  $$  \sum_{n\geq 0}\left(\sum_{\lambda\vdash n}f_\lambda^{2k} 
   \prod_{i=1}^k\prod_{u\in\lambda} (t_i+h_u^2)
 \right)\frac{x^n}{n!^{2k}} $$ 
Note that if we put each $t_i=0$ then we obtain the partition
generating function $\prod_{i\geq 1}(1-x^i)^{-1}$. The same question can
be asked with $h_u^2$ replaced with $c_u^2$ or $c_u$.
 \item Define a linear transformation $\psi\colon \Lambda_\qq
   \rightarrow \qq[t]$ by
  $$ \psi(s_\lambda) = H_\lambda^{-1}\prod_{u\in\lambda}(t+h_u^2). $$
Is there a nice description of $\psi(p_\mu)$?
\ee


\pagebreak
\begin{center}
{\large\bf Footnotes}\\[.5in]
\end{center}

Affiliation of author. Department of Mathematics, Massachusetts
Institute of Technology, Cambridge MA 02139

\medskip
$^1$2000 Mathematics Subject Classification: Primary 05E10, Secondary
05E05.  \\Key words and phrases: partition; hook length; content;
shifted part; standard Young tableau

\medskip
$^2$This material is based upon work supported by the
  National Science Foundation under Grant No.~0604423. Any opinions,
  findings and conclusions or recommendations expressed in this
  material are those of the author and do not necessarily reflect
  those of the National Science Foundation.

\end{document}